\documentstyle{article}

\def \I {{\hbox{Im }}}

\def \T {{\hbox{Tr }}}

\def \s {\sigma}

\def \g {\gamma}

\def \l {\lambda}

\def \d {{\hbox { d}}}

\def \E {{\bf E }}

\def \v {\vert}

\def \a {\alpha}

\def \de {\delta}

\def \D {\Delta}

\def \i {{\hbox {i}}}

\def \vr {\varrho}

\begin{document}
\author{A.BOUTET de MONVEL and 
A.KHORUNZHY\thanks{A.K. is grateful to the MAE (France) for the
financial support}
\\UFR de Mathematiques, Universit\'e Paris-7 Denis Diderot,
\\Paris 75251 Cedex 05, FRANCE
\\and
\\Institute for Low Temperature Physics,\\
Kharkov 310164, UKRAINE}
\title{ON UNIVERSALITY OF SMOOTHED EIGENVALUE DENSITY
OF LARGE RANDOM MATRICES}

\maketitle
\begin{abstract}
We describe the resolvent approach for the rigorous study of the
mesoscopic regime  of  Hermitian 
matrix spectra. We present results 
reflecting universal behaviour of the
smoothed density of eigenvalue distribution
of large random matrices.
\end{abstract}

Random matrices of large dimensions introduced and studied
by E.Wigner \cite{W} have applications in various
fields of theoretical physics (see e.g.  monographs
and reviews \cite{D,G,MPV} and references there).
In these studies, the spectral properties of random matrix
ensembles play an important role. 

Here the
universality conjecture
for large random matrices, formulated by F.Dyson \cite{D,D1},
is known  as the most interesting  and challenging problem.
It concerns the asymptotically local
spectral statistics, i.e. the functions that depend on certain
number
$q$ of eigenvalues of random $N\times N$ matrix $A_N$
and this number remains fixed when $N\to\infty$.

Loosely speaking, the universality conjecture states
that the local statistics regarded in the limit $N\to\infty$
do not depend on the details of the probability distribution
$P(A_N)$ of the ensemble but are determined by the symmetries
of the ensemble. For example, the expressions derived for  
local statistics of Hermitian ensembles are different from those of
real symmetric matrices.

Given Hermitian (or real symmetric) matrix $A_N$,
the distribution of its eigenvalues
$\l_1^{(N)}\le \dots\le \l_N^{(N)}$ is determined
by the normalized eigenvalue counting function
$$
\s_N(\l)\equiv \s(\l;A_N) := \#\{ \l_j^{(N)}\le \l\} N^{-1}
$$
or, equivalently, by the associated measure 
$$
\s_N(\Delta) \equiv \int_a^b \rho_N(\l) \d \l ,
\quad \Delta = (a,b)\subset {\bf R},
$$
with the formal density 
$$
\vr_N(\l) = {1\over N} \sum_{j=1}^N \de(\l-\l_j^{(N)}).
\eqno (1)
$$
The function $\s_N(\l)$ is 
called the empirical eigenvalue distribution
function.
Regarding $\s_N(\D_N)$,  it turns to be the local spectral
statistics when considered with the intervals of 
the length
$\vert\D_N\vert= O(1/N)$ as $N\to\infty$.

In general, the local spectral regime is rather hard to analyse
rigorously.
The universality conjecture is supported 
mainly for
those  ensembles of random matrices that
have explicit form of the joint  probability distribution
$\pi_N(\l_1,\dots,\l_N)$ of eigenvalues.
Starting from $\pi_N$, the same expression for $m$-point
correlation function is derived  by Dyson for the
circular ensemble of unitary random matrices (CUE) \cite{D},
by Mehta for GUE \cite{Mehta}, by Pastur
and Shcherbina for matrix models ensemble \cite{PS}.
This expression is given by the determinant of $m\times m$
matrix with the entries $\{\sin \pi (t_i-t_j)/\pi (t_i-t_j)\}$,
$i,j = 1,\dots, m$.
The same expression is derived in \cite{J}
for a random matrix ensemble with the entries
that are independent random variables, whose 
probability distribution is a
convolution
of the Gaussian distribution and the arbitrary one.

Our principal goal is to examine the presence of universality of
the spectral characteristics for those ensembles of random
matrices,  for which the explicit form of the joint eigenvalue   
distribution $\pi_N$ is unknown. 
For example,  random matrix with independent
$\pm 1$ entries falls into this class.
Our claim is that the eigenvalue density (1) smoothed over
the intervals 
 $\D_N\subset {\bf R}$ possesses
the universal properties as $N\to\infty$ provided the length
$l_N=\v\D_N\v$ satisfies
conditions $1\ll l_N\ll N$.

We determine the smoothing (or regularization) of (2) by the
formula
$$
R_N^{(\a)}(\l) := \int_{-\infty}^{\infty}
{N^\a\over 1 + N^{2\a} (\l-\l')^2} \, \vr_N(\l') \d \l'
\eqno (2)
$$
and note that in this case 
$$
R_N^{(\a)}(\l) = \I \T G_N(\l + \i N^{-\a}) N^{-1},
$$
where 
$G_N(z) = (A_N-z)^{-1}$.

According to the above definition, $\xi_N^{(\a)}(\l)$ 
with $\alpha =1$ represents 
the asymptotically local spectral statistics.
The opposite asymptotic regime when $\a=0$ is known as the global
one.
In this case the limit
$$
g(z) = \lim_{N\to\infty} \T G_N(z) N^{-1} ,
\quad \v \I z\v >0
$$
if it exists,
determines the limiting eigenvalue distribution $\s(\l)$
of the ensemble $\{ A_N\}$; that is
$$
\s(\l) = \lim_{N\to\infty} \s_N(\l), \quad 
g(z) = \int_{-\infty}^\infty (\l-z)^{-1} \d \s(\l).
$$

Regarding the global regime, 
the resolvent approach developed in papers
\cite{KP1,KP2} is proved to be
rather effective in studies of eigenvalue distribution
of large random matrices (see, for example
\cite{BK,KKP,KR}).
In this regime the limit of $g_N(z) = \T G_N(z) N^{-1}$ 
depends on the probability distribution of the ensemble,
i.e. is non-universal, as well as the 
fluctuations of $g_N(z)$ \cite{BK,KKP}.

We are interested in the behaviour of 
(2) in the 
case of $0<\a<1$. This regime 
is intermediate between the local and the global ones. It can be
called the mesoscopic regime in random matrix spectra.
For this regime, 
the  modified version of the resolvent approach 
was proposed in \cite{K} to study spectral properties of random
matrices with independent 
arbitrary distributed entries (see also \cite{BK1,BK2}).

As the further development of the resolvent approach of \cite{K}, 
we present the results concerning random matrices with
statistically dependent entries.
We consider the ensemble of random matrices
$$
H_{m,N}(x,y) = {1\over N} \sum_{\mu=1} ^m 
\xi_\mu(x) \xi_\mu (y), \quad x,y = 1,\dots,N,
\eqno (3)
$$
where the random variables $\{\xi_\mu(x)\}, x,\mu \in {\bf N}$ 
 have joint Gaussian distribution with zero mathematical
expectation  and covariance
$$
\E\{ \xi_\mu(x) \xi_\nu(y)\} = u^2 \de_{xy} \de_{\mu\nu}.
$$
Here $\de_{xy}$ denotes the Kronecker delta-symbol.
This ensemble first considered in \cite{MP} is now
of extensive use in the statistical mechanics of
disordered spin systems \cite{MPV} and in the modelling of
memory   in the theory of
neural networks \cite{HKP}.

\vskip 0.5cm
{\bf Theorem 1.} {\it Let $G_{m,N} (z) = (H_{m,N} - z)^{-1}$.
Then, for $N,m\to\infty, \, m/N\to c>0$, the random variable }
$$
R_{m,N}^{(\a)}(\l) :=  \I \T G_{m,N}(\l + \i N^{-\a}) N^{-1}
$$
{\it converges with probability 1 as 
 to the nonrandom limit
$$
\pi \vr _c(\l) = {1\over 2 \l u^2}
{\sqrt{  4cu^4 - [\l - (1+cu^2)]^2}
 }
\eqno (4)
$$
provided  $0<\a<1$ and $\l \in \Lambda_{c,u} = (u^2(1-\sqrt c)^2,
u^2 (1+\sqrt c)^2)$.}
     
\vskip 0.5cm 
{\bf Theorem 2.}
{\it Consider $k$  random variables
$$
\g_{m,N}^{(\a)}(i) := N^{1-\a} \left[ 
R_{m,N}^{(\a)}(\l_i) - \E R_{m,N}^{(\a)}(\l_i)\right],
i= 1,\dots, k,
$$
where $\l_i = \l + \tau_iN^{-\a}$ with given $\tau_i$. Then under
the conditions of Theorem 1 the joint distribution of the vector
$(\g_N(1),\dots ,\g_N(k))$
converges to the  Gaussian $k$-dimensional distribution
with zero average and covariance
$$
C(\tau_i,\tau_j) = { 4 - (\tau_i - \tau_j)^2\over
[4 + (\tau_i - \tau_j)^2]^2}.
\eqno (5)
$$
}

{\it Remark. } 
It is easy to see that if $\v\tau_1-\tau_2\v\to \infty$, then
$$
C(\tau_1,\tau_2) = -(\tau_1 - \tau_2)^{-2}(1+o(1)).
\eqno (6)
$$
This coincides with the
average value of  the Dyson's 2-point correlation function for
real symmetric matrices
considered at large distances $\v t_1-t_2\v \gg 1$ \cite{D}.

\vskip 0.2cm

To discuss these results, let us first  note  that 
theorem 1 proves existence of the smoothed density of eigenvalues
that coincides with that derived in \cite{MP} in the global regime;
$\vr_c(\l) = \s_c'(\l), \l>0$, where
$$
\s_c(\l) = \lim_{N\to\infty} \s(\l;H_{m,N}).
$$

This density obviously differs from the 
semicircle (or Wigner) distribution $\s_{w}(\l)$
$$
\s_{w}(\l) = \lim_{N\to\infty} \s(\l;W_N)
$$
where $W_N(x,y) = w(x,y)/\sqrt N$ are random 
symmetric matrices with independent identically distributed entries
with zero mathematical expectation and variance $v^2$.
This ensemble is known as the Wigner ensemble of random matrices.
It is 
known since the pioneering work of Wigner \cite{W} that
$$
\vr_{w}(\l) = \s_{w}'(\l) = {1\over 2\pi v^2}
\cases { \sqrt{ 4v^2 - \l^2} , & if $ \v \l \v \le 2v$,\cr
0, & if $ \v \l\v > 2v.$ \cr}
$$

It should be noted that in papers \cite{BK2,K}
we have proved analogues of Theorems 1 and 2 for the  Wigner
ensemble of random matrices. 
We have shown that
Theorem 2 is true  
provided $\E w(x,y)^8$ is bounded and $\a \in (0,1/8)$.
The correlation function $C(\tau_i,\tau_j)$ is given 
again by (5). 
Comparing these results, we conclude 
the fluctuations of the smoothed
eigenvalue density
do not feel  the dependence between matrix elements.

Thus, our results can be regarded as the statements
corroborating the universality conjecture for mesoscopic regime.
Namely, they show that in 
the mesoscopic regime the smoothed density of eigenvalues
$R^{(\a)}_N$ is the selfaveraging variable. It   converges as
$N\to\infty$ to the eigenvalue distribution of the ensemble and
depends on the probability distribution of the random matrix
ensemble. At the same time the fluctuations of $R^{(\a)}_N$
in the limit $N\to\infty$ coincide for two such different classes
of random matrices as (3) and the Wigner one.

This dual-type behaviour of the eigenvalue density 
of random matrices in the
mesoscopic  regime is well-known in theoretical physics (see,
for example the review \cite{G}, Chapter 8). 
For example, the universal properties of the 
mesoscopic 
eigenvalue density
 is studied in  \cite{AJM,BZ} for the matrix models ensemble.
It was shown that
the correlation function of the 
eigenvalue density (in theoretical physics terms, the "wide"
correlator)
depends on the edges of the spectrum. For the case of the 
symmetric support of the limiting eigenvalue distribution, 
the expression for the "wide" correlator 
coincides with the asymptotic expression
 (6). It should be noted that our result (6) does not depend
 on the support of $\vr_c(\l)$.

Let us describe the method developed for the proof of Theorems
1 and 2 (the full version will be published elsewhere). 
It represents 
a modification of the  resolvent approach proposed in papers
\cite{KP1,KP2}. 
This approach was developed to study the eigenvalue distribution
of random matrices and random operators in the global regime
$\v \I z \v >0$ as $N\to\infty$.
It is based on derivation and
asymptotic analysis of the system of relations for the moments
$L_k^{(N)} = \E [g_N(z)]^k, k\ge 1$. 
These relations
are of the following form 
$$
L_k = a L_{k-1} + b L_{k+1} +
\Phi_k^{(N)},
\eqno (7)
$$
where  the terms $\Phi_k^{(N)}$ 
can be estimated by $N^{-1} \v \I z\v^{-k}$.  This a priori
estimate implies that in the study of the asymptotic behaviour of
$g_N(z)$, one can restrict oneself with only two first relations.
Namely, all information about the limiting behaviour of 
$L_k^{(N)}, \ k\ge 1$ can be derived from relations for $L_1^{(N)}$
and $L_2^{(N)}$.

To consider $g_N(z)$ in the mesoscopic regime, we start with the
same system of relations for $L_k^{(N)}$. The main observation
made in \cite{K} is that in this case we need all the infinite
system of relations. More precisely, the more close $\a $ is to 1,
the greater number $K(\a)$ is such that we need to consider
relations for $L_k^{(N)}, k\ge K(\a)$. 

The matter is that 
 the terms  $\Phi_k^{(N)}$ can be estimated in terms of 
 $L_k^{(N)}$ multiplied by $N^{-\beta}$, where $\beta =
\min\{\a,1-\a\}$. 
The structure of relations (7) is such that $L_j^{(N)}, \, j< k$
enter into relation for $L_k^{(N)}$ with the factor
$N^{\beta(k-j)}$.
Therefore, admitting a priori estimate
$\v L_1^{(N)}\v\le N^\a$, we deduce that it enters into the
relation (7) with the factor $N^{-k\beta}$. Regarding, $k\ge
K(\a)$, 
 one obtains the relations with the terms that converge to finite
limits as
$N\to\infty$.

\vskip 0.5cm
{\bf Acknowledgments.} The authors are grateful for 
A.Its, P.Bleher, and H.Widom and other organizers of the semester
``Random Matrix Models and Their Applications'' at MSRI (Berkeley)
for the kind invitations to participate the workshops
and for the financial support.


\begin{thebibliography}{99}


\bibitem{AJM} J Ambj{\o}rn, J Jurkiewicz, Yu M Makeenko,
{\it Phys. Lett. B} {\bf 251} (1990) 517

\bibitem{BK} A Boutet de Monvel and A Khorunzhy,
{\it Markov Proc. Rel. Fields} {\bf 4}
(1998) 175-197

\bibitem{BK1} A Boutet de Monvel and A Khorunzhy, {\it Rand.
Oper. Stoch. Equations} {\bf  7} (1999) 1-22
\bibitem{BK2} A Boutet de Monvel and A Khorunzhy, {\it Rand.
Oper. Stoch. Equations}
{\bf  7} (1999), to appear


\bibitem{BZ} E Br\'ezin and A Zee, {\it Nucl. Phys. B } {\bf 402 }
(1993) 613-627

\bibitem{D} F Dyson, {\it J. Math. Phys.} {\bf 3} (1962) 140-175

\bibitem{D1} F Dyson; {\it J. Math. Phys.} {\bf 13} (1972) 90 



\bibitem{G} T Guhr, A M\"uller-Groeling, and H A Weidenm\"uller
{\it Phys. Rep.} {\bf 299 } (1998)

\bibitem{HKP} J A Hertz, A Krogh, and R G Palmer. 
{\it Introduction to the theory of neural computations,} 
Addison-Wesley (1991)

\bibitem{J} K Johansson, {\it preprint}, (1998)

\bibitem{K} A Khorunzhy {\it Rand. Oper. Stoch. Equations} 
{\bf 5} (1997) 147-162

\bibitem{KKP} A Khorunzhy, B Khoruzhenko, and L Pastur 
{\it J. Math. Phys.} {\bf 37} (1996) 5033-5060 

\bibitem{KP1} A Khorunzhy and L Pastur, {\it Commun. Math. Phys.}
{\bf 153} (1993) 605-646

\bibitem{KP2} A Khorunzhy and L Pastur,
in: {\it Spectral Operator Theory and Related Topics},
Adv. Soviet Math. {\bf 19} (1994) 97-127


\bibitem{KR} A Khorunzhy and G J Rodgers, {\it Rep. Math. Phys.},
{\bf 43} (1998) 297-319



\bibitem{MP} V  Marchenko and  L Pastur, 
{\it Math. USSR Sbornik} {\bf  1} (1967) 457-483

\bibitem{Mehta} M L Mehta. {\it Random Matrices.} Acad. Press,
New York (1991)



\bibitem{MPV} M M\'ezard, G Parisi, and M Virasoro.  
{\it Spin glass theory and beyond, } World Scientific, Singapur
(1987)

\bibitem{PS} L Pastur and M Shcherbina, {\it J. Stat. Phys.}
{\bf 86} (1997) 109-147





\bibitem{W} E.Wigner, {\it Ann. Math.} {\bf 62  }
(1955) 548-564


\end{thebibliography}
\end{document}